\documentclass[12pt]{article}
\usepackage{amssymb}
\usepackage{amsmath}

\newtheorem{thm}{Theorem}[section]

\newtheorem{lem}[thm]{Lemma}

\numberwithin{equation}{subsection}

\def\R{{\mathbb R}}

\def\1{{\mathbf 1}}
\def\a0{{\aleph_0}}

\def\eps{{\varepsilon}}

\newcommand{\ra}{{\rightarrow}}

\def\ne{n(\eps,B_d)}

\def\err{{\rm err}}
\def\alg{{\rm LIN}}

\title{Multivariate integration in $C^\infty([0,1]^d)$ is not strongly tractable}

\author{Jakub Onufry Wojtaszczyk
\\
Department of Mathematics, Informatics and Mechanics \\
University of Warsaw\\ ul. Banacha 2,
02-097 Warsaw, Poland\\ email: jakub.wojtaszczyk@zodiac.mimuw.edu.pl}

\date{June 9, 2003}

\begin{document}

\maketitle
\begin{abstract}
%
  It has long been known that the multivariate integration problem for
  the unit ball in $C^r([0,1]^d)$ is intractable for fixed finite~$r$.
  H.~Wo\'zniakowski has recently conjectured that this is true even if
  $r=\infty$.  This paper establishes a partial result in this direction.
  We prove that the multivariate integration problem, for infinitely
  differential functions all of whose variables are bounded by one, is
  not strongly tractable.
\end{abstract}



\section{Introduction}

Multivariate integration is a classical problem of numerical analysis that has been studied for
various normed spaces $F_d$ of functions of $d$ variables. In practical applications $d$ is
frequently very large, even in the hundreds or thousands.

Tractability and strong tractability of multivariate integration has been recently thoroughly
analyzed. These concepts are defined as follows. We consider the worst case setting and define
$\ne$ as the minimal number of function values that are needed to approximate the integral of any
$f$ from the unit ball of $F_d$ with an error threshold of $\eps$. We want to know how $\ne$
depends on~$d$ and $\eps^{-1}$.  The problem is \emph{tractable} if $\ne$ is bounded by a
polynomial in~$\eps^{-1}$ and~$d$, and \emph{strongly tractable} if the bound only depends
polynomially on~$\eps^{-1}$, with no dependence on~$d$.

The dependence of $\ne$ on $\eps^{-1}$ has been studied for many years, and bounds for $\ne$ in
terms of $\eps^{-1}$ are known for many $B_d$. For instance, let $F_d = C^r([0,1]^d)$ be the space
of $r$ times continuously differentiable functions defined on the $d$-dimensional unit cube with
the norm given as the supremum of the absolute values of all partial derivatives up to the $r$th
order. Bakhvalov proved in 1959, see \cite{book}, that there exist two positive numbers $c_{r,d}$
and $C_{r,d}$ such that $$ c_{r,d}\, \eps^{-d/r}\, \leq\, \ne\, \leq\, C_{r,d}\, \eps^{-d/r}\qquad
\forall\, \eps \in (0,1). $$ When $r$ is fixed this implies that multivariate integration in
$C^r([0,1]^d)$ is intractable. Indeed, $\ne$ cannot possibly be bounded by a polynomial in
$\eps^{-1}$ and $d$ since the exponent of $\eps^{-1}$ goes to infinity with $d$, and $\ne$ is
exponential in $d$.  However if $r$ varies, that is, when we consider the spaces
$F_d=C^{r(d)}([0,1]^d)$ with $\sup_d d/r(d) < \infty$, the behavior of $\ne$ is not known since we
do not have sharp bounds on $c_{r,d}$ and $C_{r,d}$. In fact, the best known bounds on $c_{r,d}$
are exponentially small in $d$, while the best known bounds on $C_{r,d}$ are exponentially large in
$d$, see again \cite{book}. Thus we can neither claim nor deny tractability or strong tractability
in the class $C^{r(d)}([0,1]^d)$ with $\sup_d d/r(d) < \infty$ on the basis of Bakhvalov's result.

The conjecture formulated in \cite{paper} states that multivariate integration in
$C^{r(d)}([0,1]^d)$ is intractable even if $r(d) = \infty$. That is, even when we consider
infinitely differentiable functions with all partial derivatives bounded by one, $\ne$ cannot be
bounded by a polynomial in $\eps^{-1}$ and $d$. Although we are not able to establish this
conjecture in full generality, we shall prove that multivariate integration in $C^\infty([0,1]^d)$
is not {\em strongly} tractable. This is achieved by showing that for a fixed $n$, the $n$th
minimal error goes to one as $d$ approaches infinity. More precisely, we show that for any $n$ and
$\eta$ there exists $d=d(n,\eta)$ such that for any linear algorithm there exists a polynomial that
is a sum of univariate polynomials, which belongs to the unit ball of $C^\infty([0,1]^d)$, whose
integral is at least $1-\eta$, and the algorithm outputs zero. This proof technique allows us to
prove the lack of strong tractability, but seems too weak to establish the lack of tractability.

\section{Multivariate Integration in $C^\infty([0,1]^d)$}

We precisely define the problem of multivariate integration studied in this paper. Let $F_d =
C^\infty([0,1]^d)$ be the space of real functions defined on the unit cube $[0,1]^d$ that are
infinitely differentiable with the norm $$| |f||_d\,=\,\sup\{\,|D^{\alpha}f(x)|\, :\  x \in
[0,1]^d,\  \alpha - \text{any multiindex}\,\}. $$ Here $\alpha =
[\alpha_1,\alpha_2,\ldots,\alpha_d]$ with non-negative integers $\alpha_j$, $|\alpha| = \alpha_1 +
\alpha_2 + \ldots + \alpha_d$, and $$D^\alpha f = \frac{\partial^{|\alpha|}}{\partial
x_1^{\alpha_1} \partial x_2^{\alpha_2}\ldots\partial x_d^{\alpha_d}}f$$ stands for the
differentiation operator.


Let $B_d$ denote the unit ball of this space. The multivariate integration problem is defined as an
approximation of integrals $$ I_d(f) = \int_{[0,1]^d} f(t)\, dt \qquad \forall f \in B_d, $$ where
the integral is taken with respect to the Lebesgue measure.

It is well known that adaption and the use of non--linear algorithms do not help for the
multivariate integration problem, as proven in \cite{blackbook}, see also \cite{book}. That is why
we consider only linear algorithms, $$A_{n,d}(f) = \sum_{j=1}^n a_j f(x_j)$$ for some real
coefficients $a_j$ and some sample points $x_j$ from $[0,1]^d$. Here $n$ denotes the number of
function values used by the algorithm. Of course the $a_i$ and $x_i$ may depend on $n$. The (worst
case) error of the algorithm $A_{n,d}$ is defined as $$ \err(A_{n,d}) = \sup\{|I_d(f) - A_{n,d}(f)|
:\ f \in B_d\}. $$ Let $\alg_{n,d}$ denote the class of all linear algorithms that use $n$ function
values. The $n$th minimal error is defined as $$e(n,B_d) = \inf\{\err(A_{n,d}) : \ A_{n,d} \in
\alg_{n,d}\}.$$

We shall prove the following theorem:

\begin{thm}\label{r} For any positive integer $n$ we have $$\lim_{d\ra \infty} e(n, C^\infty([0,1]^d)) = 1.$$\end{thm}


This theorem easily implies that multivariate integration in $C^\infty([0,1]^d)$ is {\em not}
strongly tractable. Indeed, were it strongly tractable, we would have a polynomial bound on $\ne$
independent of $d$, thus having a linear algorithm of error at most $\eps$ and using at most
$n=\ne$ function values. Taking, say, $\eps < \tfrac12$ and $n > n(\frac{1}{2},B_d)$ we would get
$e(n,B_d) \leq \tfrac12$ independently of $d$, and this would contradict the theorem that
$e(n,B_d)$ goes to one when $d$ approaches infinity.

\subsection{Proof of the Theorem}

We take an arbitrary positive integer $n$ and $\eta \in (0,1)$.  The idea of the proof is to
separate variables and, for sufficiently large $d$ and any $A_{n,d} \in \alg_{n,d}$, to find a
polynomial $f \in B_d$ that is a sum of univariate polynomials such that $|I_d(f) - A_{n,d}(f)| >
1-\eta$. It will suffice to find such a polynomial for which $f(x_j) = 0$ at all points $x_j$ used
by $A_{n,d}$. Then $A_{n,d}(f) = 0$, and thus $\err(A_{n,d}) \ge |I_d(f)| > 1-\eta$. Since
$A_{n,d}$ is an arbitrary linear algorithm this implies that $e(n, B_d) \geq 1 - \eta$ for
sufficiently large $d$. For the zero algorithm $A_{n,d} \equiv 0$ we have  we have $err
(A_{n,d})=1$ and hence we have $e(n, B_d) \le 1$. Since $\eta$ can be arbitrarily small, this
completes the proof of Theorem \ref{r}.

Suppose for the moment that we have the following lemma, which will be
proven in the next subsection.

\begin{lem}\label{km} For any positive integer $n$,
  and any $\eta \in (0,1)$ there exists a constant
  $K_{\eta,n}$ such that for any choice of $y_1,y_2,\ldots, y_n \in
  [0,1]$ there exists a polynomial $f : [0,1] \ra \R$ satisfying the
  following conditions:
\begin{enumerate}
\item $\max_{x\in[0,1]} |f(x)| \leq 1$,
\item $\max_{k=0,1\ldots} \max_{x \in [0,1]}|f^{(k)}(x)| \leq K_{\eta,n}$,
\item $\int_0^1 f(x)\, dx > 1-\eta$,
\item $f(y_j) = 0$ for $j = 1,2,\ldots, n$.\end{enumerate}
\end{lem}

Having Lemma \ref{km}, we take any $d \geq K_{\eta,n}$ and any $A \in
\alg_{n,d}$ that uses sample points
$$
x_j\, =\, [x_j^1, x_j^2, \ldots,x_j^d]\, \in\, [0,1]^d
$$
for $j = 1,2,\ldots, n$.

For $i = 1,2,\ldots,d$, let $f_i$ be the polynomial given by Lemma \ref{km} for $y_j = x_j^i$, with
$j=1,2,\ldots,n$. Consider the multivariate polynomial $$f(t_1,t_2,\ldots,t_d) = \frac{1}{d}
\sum_{i=1}^d f_i(t_i) \qquad t_i \in [0,1].$$ The values of $f$ are bounded by $1$ since they are
arithmetic means of the values of $f_i$ from $[-1,1]$. Any mixed derivative of such a function $f$
is 0, while $$\left|\frac{\partial^k}{\partial x_i^k}f(a)\right|\, =\, \left|\frac{1}{d}
f_i^{(k)}(a^i)\right|\, \leq\, \frac{K_{\eta,n}}{d} \, \leq\, 1. $$ Thus $f$ belongs to $B_d$.
Additionally, $$ \int_{[0,1]^d} f(t) \,dt\, =\, \frac{1}{d} \sum_{i=1}^d \int_{[0,1]} f_i(x)\, dx\,
>\, 1-\eta. $$ Furthermore $f(x_j) = 0$ since $f_i(x_j^i) = 0$ for all $j = 1,2,\ldots, n$. Thus
$f$ is a function we needed to prove Theorem \ref{r}.

\subsection{Proof of the Lemma}

We will use the Stone-Weierstrass theorem to find a function satisfying Lemma \ref{km}.

For $\delta = {\eta}/{(7n)}$, let the function $g : [0,2+\delta] \ra \R$
be defined as $1-\delta$ on $[0,1 - \delta] \cup [1 + \delta,
2+\delta]$, $-2\delta$ at 1 and linear on $[1-\delta,1]$ and
$[1,1+\delta]$. It is obviously a continuous function, so by the
Stone-Weierstrass theorem we can approximate it by a polynomial $P$ of
degree $N = N(\eta,n)$ such that $$\max_{x \in [0,2+\delta]}
|g(x)-P(x)| < \delta.$$
The polynomial $P$ is negative at 1 and
positive at $1+\delta$, so it has a root at some $y_0 \in
(1,1+\delta)$. Let $P_i(x) = P(x + y_0 - y_i)$. As $y_0 - y_i \in
(0,1+\delta)$ the polynomial $P_i$ satisfies $\max_{x \in [0,1]}
|P_i(x) - g(x +y_0-y_i)| < \delta$. Now take $$f(x) = \prod_{i=1}^n
P_i^2(x) \qquad \forall x \in [0,1].$$
Note that $f(y_j) =
\prod_{i=1}^nP^2(y_j + y_0 - y_i) = 0$ since the $j$th factor is
$P^2(y_0) = 0$.

The polynomial $P_i$ satisfies
$$
1-2\delta\,<\, P_i(x) \,<\, 1 \qquad \forall x \in [0,1]
\setminus [y_i - 2\delta,y_i + \delta].
$$
Thus
$$(1-2\delta)^{2n}\,<\,f(x)\,<\,1 \qquad \forall x \in [0,1] \setminus
\bigcup_{i=1}^n [y_i - 2\delta,y_i+\delta]$$
and, of course, $f(x) \geq
0$ on the whole interval $[0,1]$. This allows us to approximate the
integral of $f$. Indeed,
$$
\int_0^1 f(x)\, dx \,\geq\,
 \int_{[0,1] \setminus
  \bigcup_{i=1}^n [y_i - 2\delta,y_i+\delta]}\, f(x)\, dx\, >\,
(1 - 3\delta n)(1 - 2\delta)^{2n}. $$ Using the Bernoulli inequality we conclude that the last
expression is at least $$ (1 - 3n\delta) (1-4n\delta)\, \geq \,1 - 7n\delta = 1 - \eta.$$

The function $f$ is a polynomial on $[0,1]$. Its coefficients are
continuous functions of $(y_1,y_2,\ldots,y_n) \in [0,1]^n$ since the
coefficients of each $P_i$ are continuous functions of $y_i$, and $f$
is the product of $P_i$'s. The upper bound of the $j$th derivative of
$f$ is a continuous function of $f$'s coefficients, and thus a
continuous function of $(y_1,y_2,\ldots,y_n)$. As a continuous
function on a compact set it is bounded for each $j$, and so all
derivarives up to the $2nN$th order have a common bound, say,
$K_{\eta,n}$, independent of $(y_1,y_2,\ldots,y_n)$. This means that
the second condition of Lemma \ref{km} is satisfied and the proof of
Lemma \ref{km} is completed.

\subsection{Acknowledgements}
This work was completed with the enormous help of Prof.\ Henryk Wo\'{z}niakowski, who suggested the
problem, gave me the theoretical tools to solve it and worked with me on creating a publishable
paper from the mathematical reasoning. The help and support of my father was also crucial to finish
this work.

\end{document}